\title{ ~~\\ Values of the Euler phi function not divisible by
a prescribed odd prime}
\author{Pieter Moree}
\documentclass[12pt]{article}
\usepackage{amssymb, latexsym, amsfonts}
\def\@ptsize{2}
\newtheorem{Thm}{Theorem}

\newtheorem{lem}{Lemma}

\newtheorem{prop}{Proposition}
\newcommand{\qed}{\hfill $\Box$}

\begin{document}
\date{}
\maketitle
{\def\thefootnote{}
}
\begin{abstract}
\noindent 
Let $\varphi$ denote Euler's phi function. For a fixed odd
prime $q$ we give an asymptotic series expansion in the
sense of Poincar\'e for the
number ${\cal E}_q(x)$ of $n\le x$ such that $q\nmid \varphi(n)$. Thereby we
improve on a recent theorem by B.K.~Spearman and K.S.~Williams 
[Ark. Mat. {\bf 44} (2006), 166--181].
Furthermore we
resolve, under the Generalized Riemann Hypothesis, which of two approximations to 
${\cal E}_q(x)$ is asymptotically superior using recent results of Y. Ihara
on the Euler-Kronecker constant of a number field.
\end{abstract}
\section{Introduction}
Let $\varphi$ denote Euler's phi function. For a fixed odd
prime $q$ we set 
${\cal E}_q=\{n|q\nmid \varphi(n)\}$ and let ${\cal E}_q(x)$ denote the
associated counting function. (If $\cal A$ is any set of integers, then
by the associated counting function ${\cal A}(x)$ we denote the cardinality
of the elements $a$ in $\cal A$ such that $a\le x$.)
Spearman and Williams \cite{SW} proved that, as $x$ tends to infinity,
\begin{equation}
\label{uno}
{\cal E}_q(x)={e(q)x\over \log^{1/(q-1)}x}\Big(1+O_{\epsilon}({1\over \log^{1-\epsilon}x})\Big),
\end{equation}
with
\begin{equation}
\label{duo}
e(q)={(q+1)(q-1)^{q-2\over q-1}\Gamma({1\over q-1})\sin({\pi\over q-1})\over
2^{q-3\over 2(q-1)}q^{3(q-2)\over 2(q-1)}\pi^{3\over 2}(h(q)R(q)C(q))^{1\over q-1}},
\end{equation}
where $h(q)$ denotes the class number of the cyclotomic field $K(q):=\mathbb Q(\zeta_q)$ and
$R(q)$ its regulator. 
Spearman and Williams gave a rather involved description of $C(q)$, see
Section \ref{slechtec}, but we will show
that actually $C(q)=C(q,1)$, where for Re$(s)>1/2$,
\begin{equation}
\label{constant}
C(q,s)=\prod_{p\ne q\atop f_p\ge 2}\Big(1-{1\over p^{sf_p}}\Big)^{q-1\over f_p},
\end{equation}
where the sum is over all primes $p\ne q$ such that $f_p$, the smallest integer
$k\ge 1$ such that $p^k\equiv 1({\rm mod~}q)$, satisfies $f_p\ge 2$.
One has $C(3)=\prod_{p\equiv 2({\rm mod~}3)}(1-1/p^2)$ for example (this is
Lemma 3.1 of \cite{SW}).

The goal of this note is to point out that the theory of Frobenian functions allows one
to prove an estimate for ${\cal E}_q(x)$ which is much more precise than (\ref{uno}), namely 
(\ref{zozithet}). Moreover, we will show that making use of the Euler product for
the Dedekind zeta function of a cyclotomic number field, cf. (\ref{rezzie}), leads to a 
simplification of the arguments of Spearman and Williams. It allows one for example to
infer that $C(q)=C(q,1)$ and to give a very short proof of the estimate (\ref{prodigee}).

The theory
of Frobenian functions was initiated by Landau \cite{L} (and, independently, but only heuristically,
by Ramanujan \cite{MC}), continued by Bernays (of later
fame in logic) in his PhD thesis and much later by Serre \cite{S} and brought in its present state by 
Odoni, see e.g. \cite{O}.

For our purposes Theorem \ref{serre}, more or less implicit in the work of Landau already, will
do. Before stating it, we first define what a Frobenian set of primes is.
A set of primes $\cal P$ is called Frobenius of density $\delta$, if there exists a finite
Galois extension $K/\mathbb Q$ and a subset $H$ of $G:={\rm Gal}(K/\mathbb Q)$ such that
\begin{itemize}
\item $H$ is stable under conjugation;
\item $|H|/|G|=\delta$;
\item for every prime $p$, with at most finitely many exceptions, one has $p$ in $\cal P$
if $\sigma_p(K/\mathbb Q)$ is in $H$, where $\sigma_p(K/\mathbb Q)$ denotes the Frobenius
map of $p$ in $G$ (defined modulo conjugation in case $p$ does not divide the discriminant of $K$).
\end{itemize}
\begin{Thm}
\label{serre}
{\rm \cite{S}}. Let $E$ be a set of integers and $E'$ its complement in the set of natural numbers. Suppose
that $E'$ is multiplicative, that is if $a$ and $b$ are coprime positive integers, then
$$ab\in {E'} \iff \{a\in E' {\rm ~or~}b\in E'\}.$$
Put $h(s)=\sum_{n\in E'}n^{-s}$. Let $P$ be the set of primes that are in $E$.
Suppose that $P$ is Frobenian of density $\delta$,
with $0<\delta<1$. Then $h(s)/s$ has an expansion around the point $s=1$ of the form
$${h(s)\over s}={1\over (s-1)^{1-\delta}}(c_0+c_1(s-1)+\cdots+c_k(s-1)^k+\cdots).$$
Furthermore, for every integer $k\ge 2$ we have
$$E'(x)={x\over \log^{\delta} x}\Big(e_0 + {e_1\over \log x}+\cdots+{e_k\over \log^k x}+O({1\over \log^{k+1}x})\Big),$$
with $e_j=c_j/\Gamma(1-j-\delta)$.
\end{Thm}

In our problem at hand it turns out that $P$ is the set of primes $p\equiv 1({\rm mod~}q)$. But this
is precisely the set of primes $p$ that split completely in $K(q)$ and thus $\zeta_{K(q)}(s)$, 
the Dedekind zeta function of $K(q)$, comes
into play. 
We put $\alpha(q):={\rm Res}_{s=1}\zeta_{K(q)}(s)$.
The reader unfamiliar with this material is referred to Section \ref{prel}.
\begin{Thm}
\label{main}
Let $q$ be an odd prime.
Put
\begin{equation}
\label{fqs}
h_q(s)={(1-q^{-2s})\zeta(s)\over (C(q,s)(1-q^{-s})\zeta_{K(q)}(s))^{1\over q-1}}.
\end{equation}
Then $h_q(s)/s$ has an expansion around the point $s=1$ of the form
$${h_q(s)\over s}={1\over (s-1)^{(q-2)/(q-1)}}\Big(c_0(q)+c_1(q)(s-1)+\cdots+c_k(q)(s-1)^k+\cdots \Big),$$
For any $k\ge 2$ we have
\begin{equation}
\label{zozithet}
{\cal E}_q(x)={x\over \log^{1/(q-1)} x}\Big(e_0(q) + 
{e_1(q)\over \log x}+\cdots+{e_k(q)\over \log^k x}+O({1\over \log^{k+1}x})\Big),
\end{equation}
where $e_j(q)=c_j(q)/\Gamma({q-2\over q-1}-j)$.
In particular,
$$e_0(q)={(1-{1\over q^2})\over \Gamma({q-2\over q-1})(C(q)(1-{1\over q})\alpha(q))^{1\over q-1}}.$$
\end{Thm}

On using that 
$$\Gamma({1\over q-1})\Gamma({q-2\over q-1})={\pi\over \sin{\pi\over q-1}}$$
and formula (\ref{residu}), it is seen after some easy computation that
$e(q)=e_0(q)$. Thus the estimate (\ref{uno}), that is the theorem of Spearman and
Williams, is a weaker form of Theorem \ref{main}.

In the second part of the paper we deal with the problem of whether the
$${e(q)x\over \log^{1/(q-1)}x}{\rm ~-naive-~or~}e(q)\int_2^x{dt\over \log^{1/(q-1)}t}{\rm 
~-Ramanujan~type-}$$
approximation yields -asymptotically- a better approximation to ${\cal E}_q(x)$. 
For every odd prime $q$ this can be decided using Theorem \ref{formulebrij}.
Using recent results of Ihara \cite{I}, which assume
the Generalized Riemann Hypothesis (GRH) to be true, we will establish the following theorem.
\begin{Thm} 
\label{maintwo} {\rm (GRH).}
Let $q$ be
an odd prime. For $q\le 67$ the Ramanujan type approximation is asymptotically better than the
naive approximation for ${\cal E}_q(x)$, for all remaining
primes the naive 
approximation is asymptotically better.
\end{Thm}

\section{Preliminaries}
\label{prel}
For a general number field $K$ we have, for Re$(s)>1$,
$$\zeta_K(s):=\sum_{\mathfrak{a}}{1\over N{\mathfrak{a}}^{s}}=\prod_{\mathfrak{p}}{1\over 1-N{\mathfrak{p}}^{-s}}.$$
The letter $\mathfrak{a}$ will be used to denote a non-zero ideal in ${\cal O}_K$, the ring of integers of $K$,
and $\mathfrak{p}$ will be used to denote a non-zero prime ideal in ${\cal O}_K$.
This function is the {\it Dedekind zeta function} of $K$. It is known that the sum and product converge for
Re$(s)>1$, that $\zeta_K(s)$ can be analytically continued to a neighborhood of $1$ (in fact, to the whole
complex plane), and that at $s=1$ it has a simple pole. 

Let $\alpha_K$ denote the residue of the pole at
$s=1$. It is known that
\begin{equation}
\label{zoklassiek}
\alpha_K={2^{r_1}(2\pi)^{r_2}h_KR_K\over w_K\sqrt{|d_K|}},
\end{equation}
where $r_1$ is the number of real infinite primes, $r_2$ is the number of complex infinite primes, $h(K)$ is the
class number of $K$, $R(K)$ is the regulator of $K$, $w(K)$ is the number of roots of unity in $K$, and $D(K)$
is the discriminant of $K$.

Around $s=1$ we have the Laurent expansion
\begin{equation}
\label{laurent}
\zeta_K(s)={\alpha_K\over s-1}+\gamma_K+\gamma_1(K)(s-1)+\gamma_2(K)(s-1)^2+\cdots
\end{equation}
The constants $\gamma_j(\mathbb Q)$ are known as the Stieltjes constants. In particular, 
we have $\gamma_K=\gamma$, with $\gamma$ 
the Euler-Mascheroni constant. The constant ${\cal EK}_K:=\gamma_K/\alpha_K$ is called the {\it Euler-Kronecker constant} in
Ihara \cite{I} and Tsfasman \cite{T}, the reason
for this being that in the case when $K$ is imaginary quadratic the well-known
Kronecker limit formula expresses $\gamma_K$ in terms of special values of the Dedekind $\eta$ function.

\subsection{Preliminaries on cyclotomic fields}
We recall some facts from the theory of cyclotomic fields needed for our proofs. For 
a nice introduction to cyclotomic fields see \cite{W}.

The following result, see
e.g. \cite[Theorem 4.16]{N}, describes the splitting of primes in the ring of integers of
a cyclotomic field.
\begin{lem} 
\label{washington}
(cyclotomic reciprocity law). Let $K=\mathbb Q(e^{2\pi i/m})$. If the rational prime $p$ does not divide
$m$ and $f$ is the least natural number such that $p^f\equiv 1({\rm mod~}m)$, then $(p)$ (considered as
an ideal in the ring of integers of $K$) equals $\mathfrak{p}_1\cdots \mathfrak{p}_g$ with
$g=\varphi(m)/f$, all $\mathfrak{p_i}$'s distinct and of degree $f$.

However, if $p$ divides $m$, $m=p^am_1$ with $p\nmid m_1$ and $f$ is the least positive integer such that
$p^{f}\equiv 1({\rm mod~}m_1)$, then $(p)=(\mathfrak{p}_1\cdots \mathfrak{p}_g)^e$ with $e=\varphi(p^a)$, 
$g=\varphi(m_1)/f$, all $\mathfrak{p_i}$'s being distinct and of degree $f$.
\end{lem}

In case $K=K(q)$, we have $R(K)=r(q)$, $h(K)=h(q)$, $r_1=0$, $r_2=q-1$, $w(K(q))=2q$ (as $K(q)$ contains
exactly $\{\pm 1,\pm \omega_q, \pm \omega_q^2, \ldots, \pm \omega_q^{q-1}\}$ as roots of unity, with
$\omega:=e^{2\pi i/(q-1)}$)
and furthermore $D(K(q))=(-1)^{q(q-1)/2}q^{q-2}$, and thus we obtain from (\ref{zoklassiek}) that
 \begin{equation}
\label{residu}
\alpha(q)={\rm Res}_{s=1}\zeta_{K(q)}(s)=2^{q-3\over 2}q^{-{q\over 2}}\pi^{q-1\over 2}h(q)R(q).
\end{equation}

For cyclotomic fields $K(q)$ the Euler product for $\zeta_{K(q)}(s)$ can be written down explicitly using
the ``cyclotomic reciprocity law''. We find that
\begin{eqnarray}
\label{rezzie}
\zeta_{K(q)}(s)&=&(1-{1\over q^s})^{-1}\prod_{p\ne q}\Big(1-{1\over p^{sf_p}}\Big)^{1-q\over f_p}\cr
&=&(1-{1\over q^s})^{-1}\prod_{p\equiv 1({\rm mod~}q)}\Big(1-{1\over p^s}\Big)^{1-q}
\prod_{p\ne q,~p\not\equiv 1({\rm mod~}q)}\Big(1-{1\over p^{sf_p}}\Big)^{1-q\over f_p}.
\end{eqnarray}
Thus we can write
\begin{equation}
\label{gsm}
\zeta_{K(q)}(s)=(1-{1\over q^s})^{-1}{g(s)^{1-q}\over C(q,s)},
\end{equation}
where
$$g(s):=\prod_{p\equiv 1({\rm mod~}q)}\Big(1-{1\over p^s}\Big).$$

Let $k$ be a natural number and $\chi$ a character modulo $k$. Let $\chi_0$ be the principal
character modulo $k$. The Dirichlet L-series corresponding to $\chi$ is given by
$$L(s,\chi)=\sum_{n=1}^{\infty}{\chi(n)\over n^s},$$
where $s=\sigma+it\in \mathbb C$. For $\chi\ne \chi_0$ the latter series converges for $\sigma>0$
and
$$L(1,\chi)=\sum_{n=1}^{\infty}{\chi(n)\over n}=\prod_p\Big(1-{\chi(p)\over p}\Big)^{-1}\ne 0.$$
Let $g$ be a primitive root modulo $q$. For any integer not divisible by $q$ we define
the index, denoted by ind$_g(n)$, of $n$ with respect to $g$ modulo $q-1$ by
$$n\equiv g^{{\rm ind}_g(n)}({\rm mod~}q).$$
Associated with $g$ we define a character $\chi_g$ modulo $q$ by
$$\chi_g(n)=\cases{ \omega^{{\rm ind}_g(n)} & if $q\nmid n$;\cr
0 & otherwise.}$$
There are exactly $\varphi(q)=q-1$ characters modulo $q$. They are 
$\chi_0,~\chi_g,~\chi_g^2,\cdots,\chi_g^{q-2}$,
where $\chi_g^{q-1}=\chi_0$, the trivial character.

It is well-known that
\begin{equation}
\label{product}
\zeta_{K(q)}(s)=\zeta(s)L(s,\chi_g)L(s,\chi_g^2)\cdots L(s,\chi_g^{q-2}).
\end{equation}

\section{The constant $C(q)$}
\label{slechtec}
In this section, for the convenience of the reader, we repeat the definition of Spearman and
Williams of $C(q)$ and, moreover, we will show that $C(q)=C(q,1)$.

Spearman and Williams put
$$C(q,r,\chi_q):=\prod_{\chi_g(p)=\omega^r}\Big(1-{1\over p^{(q-1)/(r,q-1)}}\Big),$$
where the product is taken over all primes $p$ such that $\chi_g(p)=\omega^r$ and $(r,q-1)$ denotes
the greatest common divisor of $r$ and $q-1$.
Then they define 
$$C(q):=\prod_{r=1}^{q-2}C(q,r,\chi_g)^{(r,q-1)}.$$
From this definition it is not a priori clear that it does not depend on the choice
of the primitive root $g$. However, Spearman and Williams show that it indeed does not
depend on the choice of the primitive root $g$.
\begin{prop}
\label{cc}
We have $C(q)=C(q,1)$.
\end{prop}
{\it Proof}. By the definition of Spearman and Williams we have
\begin{equation}
\label{proddie}
C(q)=\prod_{r=1}^{q-2}\prod_{\chi_g(p)=\omega^r}
{\Big(1-{1\over p^{(q-1)/(r,q-1)}}\Big)}^{(r,q-1)}.
\end{equation}
We claim that if $\chi_g(p)=\omega^r$, then $f_p=(q-1)/(r,q-1)$. We have
$1=\chi_g(p^{f_p})=\omega^{rf_p}$. It follows that $(q-1)|rf_p$ and thus $q_r:=(q-1)/(r,q-1)$ must
be a divisor of $f_p$. On the other hand since $\chi_g(a)=1$ iff $a$ is the identity, it follows
from $\omega^{rq_r}=\chi_g(p^{q_r})=1$ and $q_r|f_p$, that $f_p=q_r$. Thus we can rewrite
(\ref{proddie}) as
\begin{equation}
\label{proddie2}
C(q)=\prod_{r=1}^{q-2}\prod_{\chi_g(p)=\omega^r}\Big(1-{1\over p^{f_p}}\Big)^{q-1\over f_p}.
\end{equation}
Note that $p\ne q$ and $f_p\ge 2$ iff $\chi_g(p)=\omega^r$ for some $1\le r\le q-2$.
This observation in combination with (\ref{proddie2}) and the absolute convergence of the
double product (\ref{proddie2}), then shows that $C(q)=C(q,1)$. \qed\\

\noindent {\tt Remark}. Proposition \ref{cc} says that $1/C(q)$ is the contribution at $s=1$
of the primes $p\ne q$, $p\not\equiv 1({\rm mod~}q)$ to the Euler product (\ref{rezzie}) of the Dedekind
zeta function of $\mathbb Q(\zeta_q)$.

\section{Proof of Theorem \ref{main}} 
\noindent {\it Proof}. 
We apply Theorem \ref{serre}.
We let $E$ be the set of natural numbers $n$ such that $q|\varphi(n)$.
Then the couting function we are after is $E'(x)$. The multiplicativity of $\varphi$
ensures that $E'$ is a multiplicative set. The set of primes $P$ in $E$ consists
of all primes $p$ such that $p\equiv 1({\rm mod~}q)$. This set is Frobenian: it
consists precisely of the primes $p$ that split completely in $K(q)$.
We have $H=$id, $G\cong (\mathbb Z/q\mathbb Z)^*$, $\delta=1/(q-1)$ and hence $0<\delta<1$.
Thus all conditions of Theorem \ref{serre} are met.

It remains to determine the $c_i$. For this we have to compute 
$h(s)=\sum_{n=1}^{\infty}f_q(n)n^{-s}$, where $f_q(n)=1$ if
$q\nmid \varphi(n)$ and $f_q(n)=0$ otherwise. Note, cf. Proposition 5.1 
of \cite{SW}, that $n\in E'$, that is $q\nmid \varphi(n)$ iff $n=\prod_{p\not\equiv 1({\rm mod~}q)}p^{e_p}$
or $n=q\prod_{p\not\equiv 1({\rm mod~}q)}p^{e_p}$, where the product is taken over all primes $p\not= q$
with $p\not\equiv 1({\rm mod~}q)$ and the $e_p$ are non-negative integers.
Thus we have, for Re$(s)>1$,
\begin{eqnarray}
\label{basisf}
h(s)&=&(1+{1\over q^s})\prod_{p\not\equiv 1({\rm mod~}q)\atop p\ne q}{1\over 1-p^{-s}}
=(1-{1\over q^{2s}})\prod_{p\not\equiv 1({\rm mod~}q)}{1\over 1-p^{-s}}\cr
&=&(1-q^{-2s})\zeta(s)g(s).
\end{eqnarray}
On using (\ref{gsm}) to express $g(s)$ in terms of $\zeta_{K(q)}(s)$, we find
that $h(s)=h_q(s)$. Now invoke Theorem \ref{serre} together with Proposition \ref{cc}. \qed

\section{On the second order coefficient}
\subsection{A comparison problem}
The prime number states that asymptotically $\pi(x)$, the number of primes $p\le x$, satisfies
$\pi(x)\sim x/\log x$. It is well-known that 
Li$(x):= \int_2^x{dt/\log t}$, the logarithmic integral, yields a much
better approximation to $\pi(x)$. Likewise one might wonder whether
$${\cal N}_q(x):={e(q)x\over \log^{1/(q-1)}x}{\rm ~or~}{\cal R}_q(x):=e(q)\int_2^x{dt\over \log^{1/(q-1)}t}$$
yields -asymptotically- a better approximation to ${\cal E}_q(x)$. The former approximation we
will call the `{\it naive approximation}' and the second the `{\it Ramanujan type approximation}' to
${\cal E}_q(x)$. To be more precise, we say that ${\cal N}_q(x)$ is a better approximation to ${\cal E}_q(x)$
than ${\cal R}_q(x)$ if 
\begin{equation}
\label{compare}
|{\cal N}_q(x)-{\cal E}_q(x)|<|{\cal R}_q(x)-{\cal E}_q(x)|,
\end{equation}
for all $x$ sufficiently large.

In the history of the theory of Frobenian functions Ramanujan was the first to put forward a problem of this
type.
If $B(x)$ denotes the counting function of integers $n\le x$ that can be written as sum
of two squares, Ramanujan in his first letter (16 Jan. 1913) to Hardy claimed that, for every $r\ge 1$,
\begin{equation}
\label{boehoe}
B(x)=K\int_2^{x}{dt\over \sqrt{\log t}}+O({x\over \log^r x}),
\end{equation}
where $K$ is a certain constant, now called the Landau-Ramanujan constant. Landau \cite{L} had
proved in 1908 that $B(x)\sim Kx/\sqrt{\log x}$: a much weaker assertion.

There is some evidence that along with his final letter (12 Jan. 1920) to Hardy, Ramanujan included
a manuscript on congruence properties of $\tau(n)$ and $p(n)$, the partition function. In this
manuscript, see \cite{BO}, Ramanujan considers, for various special primes $q$, the quantity
$\sum_{n\le x,~q\nmid \tau(n)}1$ and makes claims similar to (\ref{boehoe}). He defines $t_n=1$ if
$q\nmid \tau(n)$ and $t_n=0$ otherwise. He then typically writes: ``It is easy to prove by quite
elementary methods that $\sum_{k=1}^n t_k=o(n)$. It can be shown by transcendental methods
that
\begin{equation}
\label{boehoe1}
\sum_{k=1}^n t_k=C_q\int_2^n {dx\over (\log x)^{\delta_q}}+O\Big({n\over (\log n)^{r}}\Big),
\end{equation}
where $r$ is any positive number''. In each case he gave specific values of $C_q$ and $\delta_q$.

The above claims (\ref{boehoe}) and (\ref{boehoe1}) where shown to be asymptotically correct by Landau, respectively
Rankin. Shanks \cite{Sh} showed that (\ref{boehoe}) is false for every $r>2$. Likewise, the author \cite{M}
showed that all claims of the format (\ref{boehoe1}) in the unpublished manuscript to be false for every $r>2$.
The proof involves computing the Euler-Kronecker constant for the generating series 
$\sum_{k=1}^{\infty}t_kk^{-s}$ with several decimals of accuracy.

The comparison problem (\ref{compare}) can be studied by computing the second order coefficient $e_1(q)$
in (\ref{zozithet}) with enough precision. So that is what we set out to do. This will require an
excursion in generalized von Mangoldt functions associated to multiplicative functions. 

\subsection{Generalized von Mangoldt functions}
Let $f$ be a nonnegative real-valued multiplicative function. 
We denote the formal Dirichlet series $F(s):=\sum_{n=1}^{\infty}f(n)n^{-s}$ associated
to $f$ by $L_f(s)$. We define $\Lambda_f(n)$ by
$$-{L_f'(s)\over L_f(s)}=\sum_{n=1}^{\infty}{\Lambda_f(n)\over n^s}.$$
The notation suggests that $\Lambda_f(n)$ is an analogue of the von Mangoldt function.
Indeed, if $f={\bf 1}$, then $L_{f}(s)=\zeta(s)$ and $\Lambda_f(n)=\Lambda(n)$.
\begin{Thm}
\label{vier}
In case $f$ is a multiplicative function satisfying
$0\le f(p^r)\le c_1c_2^r$, $c_1\ge 1$, $1\le c_2<2$ and
$\sum_{p\le x}f(p)=\tau{\rm Li}(x)+O(x\log^{-2-\rho}x)$, where $\tau$ and $\rho$ are
positive real fixed numbers, then there
exists a constant $B_F$ such that
$$\sum_{n\le x}{\Lambda_f(n)\over n}=\tau \log x + B_f+O(\log^{-\rho}x).$$
Moreover, we have
$$\sum_{n\le x}f(n)=\lambda_1(f)x\log^{\tau-1}x\Big(1+(1+o(1)){\lambda_2(f)\over \log x}\Big),$$
where $\lambda_2(f)=(1-\tau)(1+B_f)$.
Alternatively we have
\begin{equation}
\label{verbinding}
B_{f}=-\lim_{s\rightarrow 1+0}\Big({L_{f}'(s)\over L_{f}(s)}+{\tau\over s-1}\Big).
\end{equation}
\end{Thm}
{\it Proof}. Identity (\ref{verbinding})  is Lemma 1 of \cite{M}. The remainder
is Theorem 4 of \cite{M2}. \qed\\

Recall that $f_q(n)=1$ if $q\nmid \varphi(n)$ and $f_q(n)=0$ otherwise.
Note that $L_{f_q}(s)=h_q(s)$. Using equation (\ref{basisf}) we find
$$-{L_{f_q}'(s)\over L_{f_q}(s)}={\log q\over q^s+1}+\sum_{p\ne q,~p\not\equiv 1({\rm mod~}q)}                     
{\log p\over p^s-1}=-{2\log q\over q^{2s}-1}+\sum_{p\not\equiv 1({\rm mod~}q)}                     
{\log p\over p^s-1},$$
from which we infer that
\begin{equation}
\label{wahn}
\Lambda_{f_q}(n)=\cases{(-1)^{r+1}\log q & if $n=q^r$, $r\ge 1$;\cr
\log p & if $n=p^r,~p\ne q,~p\not\equiv 1({\rm mod~}q)$, $r\ge 1$;\cr
0 & otherwise.}
\end{equation}
We have $h_q(s)=\sum_{n=1}^{\infty}f_q(n)n^{-s}$. By the prime number theorem for arithmetic progressions one
has 
$\sum_{p\le x}f_q(p)=\tau_q{\rm Li}(x)+O_q(x\log^{-2-\rho}x)$, with $\tau_q=(q-2)/(q-1)$ and $\rho>0$
arbitrary. We thus infer by Theorem \ref{vier} that
\begin{equation}
\label{val1}
\sum_{n\le x}{\Lambda_{f_q}(n)\over n}=\tau_q\log x +B_{f_q} + O_{\rho,q}(\log^{-\rho}x).
\end{equation}
From Theorem \ref{vier}, (\ref{zozithet}) and (\ref{val1}), we infer the
following result.
\begin{lem}
\label{fformule}
Let $q$ be an odd prime. Then
$${e_1(q)\over e_0(q)}={1\over q-1}\Big(1+B_{f_q}\Big),$$
where 
\begin{equation}
\label{nogniet}
B_{f_q}=\lim_{x\rightarrow \infty}\Big(\sum_{n\le x}{\Lambda_{f_q}(n)\over n}-
({q-2\over q-1})\log x\Big).
\end{equation}
\end{lem}
The following lemma now shows that our comparison problem can be reduced to
a comparison problem for $B_{f_q}$:
\begin{lem}
The naive approximation gives an asymptotically better approximation to ${\cal E}_q(x)$ than
the Ramanujan type approximation if $B_{f_q}<-1/2$ (that is in this case inequality
(\ref{compare}) holds for all $x$ sufficiently large). If $B_{f_q}>-1/2$ it is
the other way around.
\end{lem}
{\it Proof}. The result easily follows on noting that, as $x\rightarrow \infty$,
$${\cal E}_q(x)={e(q)x\over \log^{1/(q-1)}x}\Big(1+{1\over q-1}{(1+B_{f_q})\over \log x}+
O_q\Big({1\over \log^2 x}\Big)\Big),$$
and
$$e(q)\int_2^x{dt\over \log^{1/(q-1)}t}={e(q)x\over \log^{1/(q-1)}x}\Big(1+{1\over (q-1)\log x}
+O_q\Big({1\over \log^2 x}\Big)\Big),$$
where the first estimate is a consequence of (\ref{zozithet}) and 
Lemma \ref{fformule} and the latter follows by
partial integration. \qed\\

We will work out the limit result (\ref{nogniet}) more explicitly and then use it to
approximate $B_{f_q}$.
\begin{lem}
\label{HisJ}
Put $$H_{f_q}(x)=-{2\log q\over q^2-1}+\sum_{p\le x\atop p\not\equiv 1({\rm mod~}q)}
{\log p\over p-1}-({q-2\over q-1})\log x,$$
and  
$$J_{f_q}(x)=-\gamma-{2\log q\over q^2-1}-\sum_{p\le x\atop p\equiv 1({\rm mod~}q)}
{\log p\over p-1}+{\log x\over q-1}.$$  
Then $\lim_{q\rightarrow \infty}H_{f_q}(x)=\lim_{x\rightarrow \infty}J_{f_q}(x)=B_{f_q}$.
\end{lem} 

In Table 2 one finds computations of $J_{f_q}(10^5),J_{f_q}(10^6)$ and $J_{f_q}(10^7)$. These 
computations are very easily
implemented in MAPLE, say, and give an idea of the true value of $B_{f_q}$, but 
unfortunately cannot be used to approximate $B_{f_q}$ with a prescribed degree of accuracy. To
that end we will use (\ref{verbinding}) instead of (\ref{nogniet}) in Section \ref{deto}.

The proof of Lemma \ref{HisJ} makes use of the following result.
\begin{lem}
\label{valusha}
For every $\rho>0$ we have
We have $$\sum_{p\le x}{\log p\over p-1}=\log x -\gamma+O_{\rho}(\log^{-\rho}x).$$
\end{lem}
{\it Proof}. We have
\begin{equation}
\label{boe1}
\sum_{p\le x}{\log p\over p-1}=\sum_{n\le x}{\Lambda(n)\over n}+\sum_{p\le x\atop p^r>x}{\log p\over p^r},
\end{equation}
where the sum is over all pairs $(p,r)$ with $p^r>x$.
Now note that
\begin{equation}
\label{boe2}
\sum_{p\le x\atop p^r>x}{\log p\over p^r}=O\Big({1\over x}\sum_{p\le x^{2/3}}\log p\Big)
+O({1\over x^{4/3}}\sum_{x^{2/3}<p\le x}\log p\Big)=O({1\over x^{1/3}}),
\end{equation}
where we used the estimate $\sum_{p\le x}\log p=O(x)$. We apply Theorem \ref{vier} with $f={\bf 1}$
(and so $L_f(s)=\zeta(s)$). The constant $B_f$ was first identified by de la Vall\'ee-Poussin \cite{V},
who proved that 
$$\gamma=-\lim_{x\rightarrow \infty}\Big(\sum_{n\le x}{\Lambda(n)\over n}-\log x\Big),$$
which shows that $B_f=-\gamma$. (Alternatively, we find from (\ref{laurent}) and $\gamma(\mathbb Q)=\gamma$
by logarithmic derivation that 
\begin{equation}
\label{logzeta}
{\zeta'(s)\over \zeta(s)}=-{1\over s-1}+\gamma+O(s-1).
\end{equation}
This in combination with (\ref{verbinding}) also shows that $B_f=-\gamma$.) It then
follows by Theorem \ref{vier} and the prime number theorem that
$$\sum_{n\le x}{\Lambda(n)\over n}=\log x - \gamma + O_{\rho}(\log^{-\rho}x).$$
The result now follows on combining the latter estimate with (\ref{boe1}) and (\ref{boe2}). \qed\\

\noindent {\it Proof of Lemma} \ref{HisJ}. Follows from (\ref{val1}) and (\ref{wahn}) on invoking (\ref{boe2})
and Lemma \ref{valusha}. \qed

\subsection{Determining $B_{f_q}$ using logarithmic derivation of $f_q(s)$}
\label{deto}
In this section we use (\ref{verbinding}) to calculate $B_{f_q}$. Recall that
$L_{f_q}(s)=h_q(s)$. The problem of studying the numerical behaviour of ${\cal EK}_{K(q)}$ and the
prime sum with enough accuracy will be considered in subsequent sections.
\begin{Thm}
\label{formulebrij}
We have
\begin{equation}
\label{B1}
B_{f_q}={(3-q)\log q\over (q-1)(q^2-1)}-\gamma+{{\cal EK}_{K(q)}\over q-1}+\sum_{p\ne q\atop p\not\equiv 1({\rm mod~}q)}
{\log p\over p^{f_p}-1},
\end{equation}
or alternatively,
\begin{equation}
\label{B2}
B_{f_q}={(3-q)\log q\over (q-1)(q^2-1)}-({q-2\over q-1})\gamma+{1\over q-1}\sum_{k=1}^{q-2}{L'(1,\chi_g^k)\over L(1,\chi_g^k)}
+\sum_{p\ne q\atop p\not\equiv 1({\rm mod~}q)}
{\log p\over p^{f_p}-1}.
\end{equation}
\end{Thm}
{\it Proof}. On noting that
$${1\over q-1}{C'(q,s)\over C(q,s)}=\sum_{p\ne q\atop p\not\equiv 1({\rm mod~}q)}{\log p\over p^{sf_p}-1},$$
we find by logarithmic differentation of (\ref{fqs}) that
\begin{equation}
\label{logfqs}
-{h_q'(s)\over h_q(s)}=-{2\log q\over q^{2s}-1}-{\zeta'(s)\over \zeta(s)}+
{1\over q-1}\Big({\zeta_K'(s)\over \zeta_K(s)}+{\log q\over q^s-1}\Big)
+\sum_{p\ne q\atop p\not\equiv 1({\rm mod~}q)}{\log p\over p^{sf_p}-1}.
\end{equation}
For notational convenience we put
$$v_q(s):=\sum_{p\ne q\atop p\not\equiv 1({\rm mod~}q)}{\log p\over p^{sf_p}-1}.$$
By logarithmic differentiation of the Laurent series (\ref{laurent}) we find that
\begin{equation}
\label{logderiv}
{\zeta_{K(q)}'(s)\over \zeta_{K(q)}(s)}=-{1\over s-1}+{\cal EK}_{K(q)}
+O_q(s-1).
\end{equation}
On inserting the latter estimate for $\zeta_{K(q)}'(s)/\zeta_{K(q)}(s)$ and 
the estimate (\ref{logzeta}) for $\zeta'(s)/\zeta(s)$
in (\ref{logfqs}) and
setting $\tau_q:=(q-2)/(q-1)$, we find that
$$-{h_q'(s)\over h_q(s)}-{\tau_q\over s-1}=-{2\log q\over q^{2s}-1}+{\log q\over (q-1)(q^s-1)}
+v_q(s)-\gamma+{{\cal{EK}}_{K(q)}\over q-1}+O_q(s-1).$$
Using (\ref{verbinding}) we conclude that (\ref{B1}) holds true.
By using (\ref{product}) and (\ref{logzeta}) we find similarly that (\ref{B2}) holds. (Alternatively
one merely combines (\ref{B1}) with (\ref{noidea}) to arrive at the same conclusion.)
This completes the proof. \qed\\

\noindent {\tt Remark}. Using (\ref{nul}) one can also express $B_{f_q}$ in terms of the
non-trivial zeros of $\zeta_{K(q)}(s)$.

\subsection{On the Euler-Kronecker constant for $\mathbb Q(\zeta_q)$}
Put
$${\rm Cyc}_q(x)=\log x-(q-1)\sum_{p\le x\atop p\ne q}{\log p\over p^{f_p}-1}-{\log q\over q-1}.$$
We leave it as an exercise to the reader to show that for any $\rho>0$ we have
$${\rm Cyc}_q(x)={\cal EK}_{K(q)}+O_{q,\rho}(\log^{-\rho}x).$$
By logarithmic differentation from (\ref{product}) we find that
\begin{equation}
\label{noidea}
{\cal EK}_{K(q)}=\gamma+\sum_{k=1}^{q-2}{L'(1,\chi_g^k)\over L(1,\chi_g^k)}.
\end{equation}
For example, cf. \cite{M2},
$${\cal EK}_{K(3)}=\gamma+{L'(1,\chi_{-3})\over L(1,\chi_{-3})}=0.945497280871680703239749994158189073\cdots$$
Ihara \cite{I2} conjectures that for any $\epsilon>0$ we have
$$({1\over 2}-\epsilon)\log q < {\cal EK}_{K(q)}<({3\over 2}+\epsilon)\log q,$$ for all $q$ sufficiently
large. In \cite{I} he remarks that it seems very likely that always ${\cal EK}_{K(q)}>0$. (This
was checked numerically for $q\le 8000$ by Mahora Shimura, assuming
GRH.) How large Euler-Kronecker constants
can get seems to be a much easier problem then how small they can get. 
Ihara
observed that ${\cal EK}_K$ can be conspicuously negative and that this occurs when $K$ has many primes having
small norm. However, in the case of $\mathbb Q(\zeta_q)$ the smallest norm is $q$ and therefore rather
large as $q$ increases.

Ihara has given bounds for ${\cal EK}_K$ that are valid under GRH. Specalising his result
to the case where $K=K(q)$ one obtains, on invoking the cyclotomic reciprocity law, the
following proposition.
\begin{prop} 
\label{grensje} (Ihara {\rm \cite[Proposition 2]{I}}.) Assume that GRH holds. Then we have
$
{\rm low}_q(x)\le {\cal EK}_{K(q)}\le {\rm upp}_q(x),
$
with
$$
{\rm upp}_q(x)={\sqrt{x}+1\over \sqrt{x}-1}(\log x-\Phi_{K(q)}(x)+l_q(x))+{2\kappa_q\over \sqrt{x}-1}-1,
$$
$$
{\rm low}_q(x)={\sqrt{x}-1\over \sqrt{x}+1}(\log x-\Phi_{K(q)}(x)+l_q(x))-{2\kappa_q\over \sqrt{x}+1}-1,
$$
where
$$(x-1)\Phi_{K(q)}(x)=(q-1)\sum_{p^{f_pk}\le x\atop p\ne q}\Big({x\over p^{f_pk}}-1\Big)\log p
+\sum_{q^{k}\le x}\Big({x\over q^{k}}-1\Big)\log q,~(x>1),$$
where the first sum is over all prime powers $p^{f_pk}$, $k\ge 1$, such that $p^{f_pk}\le x$ and $p\ne q$,
$$l_q(x)={(q-1)\over 2}\Big(\log{x\over x-1}+{\log x\over x-1}\Big),$$
and
$2\kappa_q=(q-2)\log q-(q-1)(\gamma+\log 2\pi)$.

Moreover, both upp$_q(x)$ and low$_q(x)$ tend to ${\cal EK}_{K(q)}$ as $x$ gets large and thus these bounds
allow one to calculate ${\cal EK}_{K(q)}$ with arbitrary precision.
\end{prop}
Theorem 1 of Ihara \cite{I} implies that, under GRH,
$${\cal EK}_{K(q)}\le \Big({z_q+1\over z_q-1}\Big)(2\log z_q+1),$$
with $z_q={(q-2)\over 2}\log q$. A simple analysis shows that
this implies that, for $q\ge 23$,
\begin{equation}
\label{ihabound}
{\cal EK}_{K(q)}\le 2\log (q\log q).
\end{equation}

As concerns the zeros $\rho$ of $\zeta_{K(q)}(s)$ in the critical strip, the 
so called non-trivial zeros, we have by \cite{veel},
\begin{equation}
\label{nul}
\sum_{\zeta_{K(q)}(\rho)=0}{1\over \rho}={\cal EK}_{K(q)}-(q-1)(\log 2 +\gamma)+
{1\over 2}(q-2)\log q -{(q-1)\over 2}\log \pi,
\end{equation}
where the zeros are counted with possible multiplicity. Since -at least conjecturally -
${\cal EK}_{K(q)}$ is small in comparison with $q\log q$ it seems to `measure' a subtle
effect in the distribution of the zeros.

In Table 1 some data concerning ${\cal EK}_{K(q)}$ are gathered. The second
and third column 
give ${\rm Cyc}_q(10^5)$, respectively ${\rm Cyc}_q(10^6)$. The next two columns
give ${\rm low}_q(x)$ and ${\rm upp}_q(x)$ for the value of $x$ recorded in the sixth
column. The final value gives the true value of 
${\cal EK}_{K(q)}$ as computed with MAGMA (computations
of $L$ and $L'$
were implemented in MAGMA by Tim
Dokchitser).

\subsection{Estimating the prime sum in $B_{f_q}$}
Put $$v(q):=\sum_{p\ne q\atop p\not\equiv 1({\rm mod~}q)}
{\log p\over p^{f_p}-1}.$$
Note that $v(q)$ is the prime sum arising in our expressions for $B_{f_q}$ in
Theorem \ref{formulebrij}. We can estimate this quantity as follows.
\begin{lem}
\label{sommieb}
For $q\ge 67$ we have
$$v(q)\le {2[\log(9\log q)][\log q]\over 3q}.$$
\end{lem}
Our proof of this makes use of the following estimate.
\begin{lem}
\label{moeizaam}
For $x\ge 3$ one has
$$\sum_{p>x}{\log p\over p^2-1}\le {1.055\over x}.$$
\end{lem}
{\it Proof}. For $x\ge 7481$ one has $0.98x\le \sum_{p\le x}\log p\le 1.017x$, as was
shown by Rosser and Schoenfeld \cite{RS}. From this one easily infers that for $x\ge 7481$
$$\sum_{p>x}{\log p\over p^k-1}\le {x\over x^k-1}\Big(-0.98+1.017{k\over k-1}\Big).$$
A simple further numerical analysis using the latter estimate with $k=2$ then gives the result. \qed\\

\noindent {\it Proof of Lemma} \ref{sommieb}. Write
$$v(q)=\sum_{p < q\atop p\not\equiv 1({\rm mod~}q)}
{\log p\over p^{f_p}-1}+\sum_{p>q\atop p\not\equiv 1({\rm mod~}q)}
{\log p\over p^{f_p}-1}=v_1(q)+v_2(q),$$
say. Note that $f_p\ge 3$ in the former sum.
If $p$ is odd, then $p^{f_p}=1+2kq$ for some integer $k$. This
together with the observation that $\log t/(t^r-1)$ is nonincreasing for $t\ge 2$
and $r\ge 2$ fixed and for $r\ge 2$ and $t\ge 2$ fixed, shows
that
$$v_1(q)\le {\log 2\over q}+\sum_{k=1}^m{1\over 3}{\log(1+2kq)\over 2kq},$$
where $m$ is the number of odd primes not exceeding $q$ (and thus $m=\pi(q)-2$). 

Now let us refine this estimate for $v_1(q)$ further. Let $g$ be the smallest integer such
that $2^g\ge q^2$. We consider the primes $p<q$
for which $f_p\le g$ first. Note that for fixed $f$ there
are at most $f-1$ primes $p<q$ with $p\not\equiv 1({\rm mod~}q)$
such that $p^f\equiv 1({\rm mod~}q)$. Note that $\sum_{f=2}^g (f-1)\le g^2/2$.
It follows that the primes $p$ with $f_p\le g$ contribute
at most
$${\log 2\over q}+\sum_{k=1}^{m_1}{1\over 3}{\log(1+2kq)\over 2kq}$$
to $v_1(q)$, where
$m_1={\rm min}\{m,[g^2/2]\}$. 
Since trivially $\pi(q)\le (q+1)/2$, we have for $q\ge 3$ that
$$1+2kq\le 1+2m_1q\le 1+2(\pi(q)-2)q<q^2$$ and hence, on using that
$\sum_{k\le n}1/k\le \log n +1$,
$$\sum_{k=1}^{m_1}{1\over 3}{\log(1+2kq)\over 2kq}\le {\log q\over 3q}\sum_{k=1}^{[g^2/2]}{1\over k}
\le {\log q\over 3q}(\log({g^2\over 2})+1)\le {\log q\over 3q}(\log({2\log ^2 q\over \log^2 2})+1).$$

The primes $2<p<q$ with $f_p>g$ contribute at most
$${1\over q^2}\sum_{p<q}\log p\le {1.0012\over q},$$ to $v_1(q)$,
where we used the estimate $\sum_{p\le x}\log p<1.0012x$ valid for $x>0$ \cite[Theorem 6]{RS2}.
We thus find that for $q\ge 3$
$$v_1(q)\le {\log 2\over q}+{\log q\over 3q}(\log({2\log ^2 q\over \log^2 2})+1)+{1.0012\over q}.$$
On invoking Lemma \ref{moeizaam} to estimate $v_2(q)$ we then find that
$$v(q)\le {\log 2\over q}+{\log q\over 3q}(\log({2\log ^2 q\over \log^2 2})+1)+{1.0012\over q}+{1.055\over
q}.$$
On some further analysis the result is easily obtained. \qed\\

\noindent {\tt Remark}. Note that in case $q$ is a Mersenne prime we have
$$v(q)\ge {\log 2\over 2^{f_2}-1}={\log 2\over q}.$$ Actually, the only
$q$ I have been able to find for which $v(q)>(\log 2)/q$
are the Mersenne primes. It thus is conceivable that if $q$ is not a 
Mersenne prime, then always $v(q)< (\log 2)/q$. For a given $\epsilon>0$
it also seems that the primes $q$ for which $v(q)>\epsilon/q$ have density zero.
In general $v(q)$ is relatively large if $q$ almost equals a number
of the form $p^r-1$ with $p$ small. For example, if $2q=3^r-1$ for some $r$ (e.g. 
when $r=3,7,13,71$),
then $v(g)>(\log 3)/(2q)$.

\subsection{Some numerical data regarding $B_{f_q}$}
It is not difficult to relate $B_{f_3}$ to the constant $B_{g_{3,2}}$ computed
with high decimal accuracy in Moree \cite[p. 437]{M2}. One finds that
$$B_{f_3}=B_{g_{3,2}}+{\log 3\over 4}=-{\gamma\over 2}+{L'(1,\chi_{-3})\over 2L(1,\chi_{-3})}
+v(3)=0.24718078879811624702914196\cdots $$

In Table 2 one finds further values of $B_{f_q}$ with 5 decimal precision.
They were computed from (\ref{B1}) using 
a precise enough approximation of ${\cal EK}_{K(q)}$
and $v(q)$. The latter was obtained using that, for $y\ge 3$,
$$\sum_{p\le y,~p\ne q\atop p\not\equiv 1({\rm mod~}q)}{\log p\over p^{f_p}-1}
< v(q)\le \sum_{p\le y,~p\ne q\atop p\not\equiv 1({\rm mod~}q)}{\log p\over p^{f_p}-1}+
{1.055\over y},$$
and taking $y$ large enough.

\subsection{The proof of Theorem \ref{maintwo}}
With the following lemma at our
disposal we are finally in the position to establish Theorem \ref{maintwo}. 
\begin{lem} {\rm (GRH)}. 
\label{katush}
Let $y\ge 3$. We have ${\rm Low}_q(x,y)\le B_{f_q}\le {\rm Upp}_q(x,y)$ with
$${\rm Low}_q(x,y)={(3-q)\log q\over (q-1)(q^2-1)}-\gamma+{{\rm low}_q(x)\over q-1}
+\sum_{p\le y,~p\ne q\atop p\not\equiv 1({\rm mod~}q)}{\log p\over p^{f_p}-1},$$
and 
$${\rm Upp}_q(x,y)={(3-q)\log q\over (q-1)(q^2-1)}-\gamma+{{\rm upp}_q(x)\over q-1}
+\sum_{p\le y,~p\ne q\atop p\not\equiv 1({\rm mod~}q)}{\log p\over p^{f_p}-1}+{1.055\over y}.$$
These estimates allow one to estimate $B_{f_q}$ with arbitrary precision.
\end{lem}
{\it Proof}. Follows from equality (\ref{B1}) on invoking Proposition \ref{grensje}
and Lemma \ref{moeizaam}. \qed\\

\noindent {\it Proof of Theorem} \ref{maintwo}. Assume GRH. On invoking (\ref{B1}) together
with the estimates (\ref{ihabound}), respectively Lemma \ref{sommieb} for ${\cal EK}_{K(q)}$ and $v(q)$, we find that
$$B_{f_q}\le -\gamma + {2 \log(q\log q)\over q}+
{2[\log(9\log q)\log q]\over 3q},$$
when $q\ge 67$.
It is easily proved that the right hand side is monotonically decreasing for $q\ge 67$. On taking $q=419$
the right hand side equals
$-0.50143\ldots$ and so we obtain that $B_{f_q}< -1/2$ for all $q\ge 419$. For $y\ge 3$ we can similarly
estimate $B_{f_q}$ as follows:
$$B_{f_q}\le -\gamma+{2\log(q\log q)\over q}+\sum_{p<y,~p\ne q\atop p\not\equiv 1({\rm mod~}q)}
{\log p\over p^{f_p}-1}+{1.055\over y}.$$
Using this estimate with $y$ chosen large enough ($y=1373$ will do) one concludes that
also $B_{f_q}<-1/2$ in the range $179\le q<419$.

For every prime $q$ in the range $71\le q\le 173$ one searches for appropriate $x$ and $y$ such that
Upp$_q(x,y)<-1/2$. By Lemma \ref{katush} it then follows that $B_{f_q}\le {\rm Upp}_q(x,y)<-1/2$
for every prime $q$ in this range.

Finally, for every prime $q\le 67$ one searches for appropriate $x$ and $y$ such that
Low$_q(x,y)>-1/2$. By Lemma \ref{katush} it then follows that $B_{f_q}\ge {\rm Low}_q(x,y)>-1/2$
for every prime $q$ in this range.\qed

\section{On Mertens' theorem for arithmetic progressions}
\noindent A crucial ingredient in the paper of Spearman and Williams is an asymptotic estimate for
$\prod_{p\le x\atop p\equiv 1({\rm mod~}q)}(1-1/p)$ \cite[Proposition 6.3]{SW}. They prove
that
\begin{equation}
\label{prodigee}
\prod_{p\le x\atop p\equiv 1({\rm mod~}q)}(1-{1\over p})
=\Big({qe^{-\gamma}\over (q-1)\alpha(q)C(q)\log x}\Big)^{1\over q-1}
\Big(1+O_q({1\over \log x})\Big).
\end{equation}
An alternative, much shorter proof of the estimate (\ref{prodigee}) can be obtained
on invoking Mertens' theorem for algebraic number fields.
\begin{lem}
\label{rosie}
Let $\alpha_K$ denote the residue of $\zeta_K(s)$ at $s=1$. Then,
$$\prod_{N\mathfrak p\le x}\Big(1-{1\over N\mathfrak p}\Big)={e^{-\gamma}\over \alpha_K\log x}
\Big(1+O_K({1\over \log x})\Big),$$
where the product is over all prime ideals $\mathfrak p$ in the ring of integers of $K$
whose norm is less than or equal to $x$.
\end{lem}
{\it Proof}. Similar to that of the usual Mertens' theorem, see e.g. Rosen \cite{R}. \qed\\

\noindent We invoke the latter result with $K=K(q)$ 
and work out the product over the prime ideals more explicitly using the cyclotomic reciprocity law, Lemma 
\ref{washington}. One finds, for $x\ge q$, that it equals
$$
(1-{1\over q})\prod_{p\le x\atop p\equiv 1({\rm mod~}q)}\Big(1-{1\over p}\Big)^{q-1}
\prod_{p\le x,~p\neq q\atop p\not \equiv 1({\rm mod~}q)}\Big(1-{1\over p^{f_p}}\Big)^{q-1\over f_p}=$$
$$(1+O_q({1\over x}))(1-{1\over q})C(q)\prod_{p\le x\atop p\equiv 1({\rm mod~}q)}\Big(1-{1\over p}\Big)^{q-1},$$
where we used that for $k\ge 2$,
$$\sum_{p>x}\Big(1-{1\over p^k}\Big)^{-1}=1+O\Big(\sum_{n>x}{n^{-k}}\Big)=1+O(x^{1-k}).$$ Thus, on 
invoking Lemma \ref{rosie} we find
$$(1-{1\over q})C(q)\prod_{p\le x\atop p\not \equiv 1({\rm mod~}q)}\Big(1-{1\over p}\Big)^{q-1}=
{e^{-\gamma}\over \alpha(q)\log x}
\Big(1+O_q({1\over \log x})\Big),$$
from which (\ref{prodigee}) is easily deduced.\\

\vfil\eject
\centerline{{\bf Table 1:} Approximate numerical values of ${\cal EK}_{K(q)}$ with $q$ an odd prime}
\begin{center}
\begin{tabular}{|c|c|c|c|c|c|c|}\hline
$q$&${\rm Cyc}_q(10^5)$&${\rm Cyc}_q(10^6)$&${\rm low}_q(x)$&${\rm upp}_q(x)$& $x$ &{\rm true}\\ \hline\hline
$3$&$0.9372\cdots$&$0.9431\cdots$&$\ge 0.945$&$\le 0.946$&$3\cdot 10^5$ & $0.94549\cdots$ \\ \hline
$5$&$1.7148\cdots$&$1.7181\cdots$&$\ge 1.719$&$\le 1.722$&$3\cdot 10^5$ & $1.72062\cdots$\\ \hline
$7$&$2.0799\cdots$&$2.0865\cdots$&$\ge 2.086$&$\le 2.090$&$10^6$ & $2.08759\cdots$ \\ \hline
$11$&$2.4216\cdots$&$2.4116\cdots$&$\ge 2.411$&$\le 2.420$&$10^6$ & $2.41542\cdots$\\ \hline
$13$&$2.6022\cdots$&$2.6050\cdots$&$\ge 2.601$&$\le 2.615$&$10^6$ & $2.61075\cdots$\\ \hline
$17$&$3.5662\cdots$&$3.5832\cdots$&$\ge 3.565$&$\le 3.592$&$10^6$ & $3.58197\cdots$\\ \hline
$19$&$4.7659\cdots$&$4.7876\cdots$&$\ge 4.765$&$\le 4.802$&$10^6$ & $4.79040\cdots$\\ \hline
$23$&$2.6185\cdots$&$2.6090\cdots$&$\ge 2.594$&$\le 2.635$&$10^6$ & $2.61128\cdots$\\ \hline
$29$&$3.0870\cdots$&$3.0932\cdots$&$\ge 3.068$&$\le 3.132$&$10^6$ & $3.09373\cdots$\\ \hline
$31$&$4.2759\cdots$&$4.3078\cdots$&$\ge 4.264$&$\le 4.340$&$10^6$ & $4.31444\cdots$\\ \hline
$37$&$4.3149\cdots$&$4.3155\cdots$&$\ge 4.262$&$\le 4.363$&$10^6$ & $4.30493\cdots$\\ \hline
$41$&$3.9661\cdots$&$3.9649\cdots$&$\ge 3.902$&$\le 4.020$&$10^6$ & $3.97152\cdots$\\ \hline
$43$&$4.3408\cdots$&$4.3802\cdots$&$\ge 4.318$&$\le 4.446$&$10^6$ & $4.37862\cdots$\\ \hline
$47$&$4.8142\cdots$&$4.7925\cdots$&$\ge 4.717$&$\le 4.865$&$10^6$ & $4.79939\cdots$\\ \hline
$53$&$4.3029\cdots$&$4.3370\cdots$&$\ge 4.267$&$\le 4.392$&$2\cdot 10^6$ & $4.33773\cdots$\\ \hline
$59$&$5.4275\cdots$&$5.4285\cdots$&$\ge 5.351$&$\le 5.501$&$2\cdot 10^6$ & $5.43351\cdots$\\ \hline
$61$&$5.0024\cdots$&$5.0618\cdots$&$\ge 4.971$&$\le 5.127$&$2\cdot 10^6$ & $5.07108\cdots$\\ \hline
$67$&$5.3340\cdots$&$5.2876\cdots$&$\ge 5.204$&$\le 5.384$&$2\cdot 10^6$ & $5.29213\cdots$\\ \hline \hline
$71$&$5.2392\cdots$&$5.2336\cdots$&$\ge 5.148$&$\le 5.343$&$2\cdot 10^6$ & $5.25525\cdots$\\ \hline
$73$&$3.9935\cdots$&$4.0650\cdots$&$\ge 3.957$&$\le 4.157$&$2\cdot 10^6$ & $4.06694\cdots$\\ \hline
$79$&$5.0581\cdots$&$5.0004\cdots$&$\ge 4.905$&$\le 5.132$&$2\cdot 10^6$ & $4.99827\cdots$\\ \hline
$83$&$2.9654\cdots$&$3.0295\cdots$&$\ge 2.900$&$\le 3.139$&$2\cdot 10^6$ & $3.03313\cdots$\\ \hline
$89$&$4.1811\cdots$&$4.1574\cdots$&$\ge 3.963$&$\le 4.341$&$10^6$ & $4.16409\cdots$\\ \hline
$97$&$4.8455\cdots$&$4.8793\cdots$&$\ge 4.660$&$\le 5.090$&$10^6$ & $4.89124\cdots$\\ \hline
$101$&$5.2782\cdots$&$5.2883\cdots$&$\ge 5.073$&$\le 5.530$&$10^6$ & $5.29701\cdots$\\ \hline
$103$&$5.1005\cdots$&$5.1326\cdots$&$\ge 4.899$&$\le 5.368$&$10^6$ & $5.14433\cdots$\\ \hline
$107$&$5.4382\cdots$&$5.5044\cdots$&$\ge 5.232$&$\le 5.728$&$10^6$ & $5.45827\cdots$\\ \hline
$109$&$6.9373\cdots$&$6.9267\cdots$&$\ge 6.664$&$\le 7.179$&$10^6$ & $6.90663\cdots$\\ \hline
$113$&$3.9793\cdots$&$4.0425\cdots$&$\ge 3.759$&$\le 4.288$&$10^6$ & $4.02173\cdots$\\ \hline
$127$&$5.0040\cdots$&$5.0705\cdots$&$\ge 4.763$&$\le 5.390$&$10^6$ & $5.08859\cdots$\\ \hline
$131$&$2.8372\cdots$&$2.8495\cdots$&$\ge 2.550$&$\le 3.917$&$10^6$ & $2.83682\cdots$\\ \hline
$137$&$4.9312\cdots$&$4.9205\cdots$&$\ge 4.607$&$\le 5.303$&$10^6$ & $4.93700\cdots$\\ \hline
$139$&$5.8719\cdots$&$5.8953\cdots$&$\ge 5.546$&$\le 6.260$&$10^6$ & $5.88916\cdots$\\ \hline
$149$&$6.0227\cdots$&$5.9895\cdots$&$\ge 5.611$&$\le 6.396$&$10^6$ & $5.98342\cdots$\\ \hline
$151$&$5.1040\cdots$&$5.0604\cdots$&$\ge 4.679$&$\le 5.474$&$10^6$ & $5.04201\cdots$\\ \hline
$157$&$7.4201\cdots$&$7.4053\cdots$&$\ge 7.007$&$\le 7.855$&$10^6$ & $7.40802\cdots$\\ \hline
$163$&$5.9314\cdots$&$5.9475\cdots$&$\ge 5.522$&$\le 6.409$&$10^6$ & $5.92966\cdots$\\ \hline
$167$&$8.1704\cdots$&$8.0129\cdots$&$\ge 7.596$&$\le 8.520$&$10^6$ & $8.03300\cdots$\\ \hline
$173$&$3.4172\cdots$&$3.3853\cdots$&$\ge 2.924$&$\le 3.874$&$10^6$ & $3.38434\cdots$\\ \hline
\end{tabular}
\end{center}

\vfil\eject
\centerline{{\bf Table 2:} Numerical data related to the evaluation of $B_{f_q}$}
\begin{center}
\begin{tabular}{|c|c|c|c|c|c|c|}\hline
$q$&$J_{f_q}(10^5)$&$J_{f_q}(10^6)$&$J_{f_q}(10^7)$&{\rm true}&$v(q)$\\ \hline\hline
$3$&$+0.2430\cdots$&$+0.2460\cdots$&$+0.2469\cdots$&$+0.24718\cdots$&$0.35164\cdots$\\ \hline
$5$&$-0.1042\cdots$&$-0.1034\cdots$&$-0.1029\cdots$&$-0.10281\cdots$&$0.07777\cdots$\\ \hline
$7$&$-0.1347\cdots$&$-0.1336\cdots$&$-0.1334\cdots$&$-0.13348\cdots$&$0.12282\cdots$\\ \hline
$11$&$-0.3419\cdots$&$-0.3429\cdots$&$-0.3425\cdots$&$-0.34255\cdots$&$0.00910\cdots$\\ \hline
$13$&$-0.3268\cdots$&$-0.3266\cdots$&$-0.3262\cdots$&$-0.32617\cdots$&$0.04620\cdots$\\ \hline
$17$&$-0.3584\cdots$&$-0.3574\cdots$&$-0.3576\cdots$&$-0.35751\cdots$&$0.00443\cdots$\\ \hline
$19$&$-0.3087\cdots$&$-0.3074\cdots$&$-0.3074\cdots$&$-0.30734\cdots$&$0.01100\cdots$\\ \hline
$23$&$-0.4627\cdots$&$-0.4631\cdots$&$-0.4630\cdots$&$-0.46308\cdots$&$0.00082\cdots$\\ \hline
$29$&$-0.4703\cdots$&$-0.4701\cdots$&$-0.4701\cdots$&$-0.47009\cdots$&$0.00034\cdots$\\ \hline
$31$&$-0.4014\cdots$&$-0.4003\cdots$&$-0.4002\cdots$&$-0.40015\cdots$&$0.03658\cdots$\\ \hline
$37$&$-0.4589\cdots$&$-0.4589\cdots$&$-0.4591\cdots$&$-0.45919\cdots$&$0.00092\cdots$\\ \hline
$41$&$-0.4797\cdots$&$-0.4797\cdots$&$-0.4794\cdots$&$-0.47957\cdots$&$0.00044\cdots$\\ \hline
$43$&$-0.4755\cdots$&$-0.4746\cdots$&$-0.4747\cdots$&$-0.47468\cdots$&$0.00021\cdots$\\ \hline
$47$&$-0.4740\cdots$&$-0.4745\cdots$&$-0.4744\cdots$&$-0.47441\cdots$&$0.00012\cdots$\\ \hline
$53$&$-0.4956\cdots$&$-0.4949\cdots$&$-0.4949\cdots$&$-0.49494\cdots$&$0.00021\cdots$\\ \hline
$59$&$-0.4847\cdots$&$-0.4846\cdots$&$-0.4845\cdots$&$-0.48460\cdots$&$0.00006\cdots$\\ \hline
$61$&$-0.4934\cdots$&$-0.4924\cdots$&$-0.4922\cdots$&$-0.49232\cdots$&$0.00143\cdots$\\ \hline 
$67$&$-0.4970\cdots$&$-0.4977\cdots$&$-0.4976\cdots$&$-0.49767\cdots$&$0.00026\cdots$\\ \hline \hline
$71$&$-0.5025\cdots$&$-0.5026\cdots$&$-0.5023\cdots$&$-0.50234\cdots$&$0.00061\cdots$\\ \hline
$73$&$-0.5211\cdots$&$-0.5201\cdots$&$-0.5200\cdots$&$-0.52013\cdots$&$0.00137\cdots$\\ \hline
$79$&$-0.5125\cdots$&$-0.5132\cdots$&$-0.5132\cdots$&$-0.51332\cdots$&$0.00049\cdots$\\ \hline
$83$&$-0.5416\cdots$&$-0.5408\cdots$&$-0.5408\cdots$&$-0.54077\cdots$&$0.00007\cdots$\\ \hline
$89$&$-0.5299\cdots$&$-0.5301\cdots$&$-0.5301\cdots$&$-0.53010\cdots$&$0.00034\cdots$\\ \hline
$97$&$-0.5270\cdots$&$-0.5266\cdots$&$-0.5265\cdots$&$-0.52657\cdots$&$0.00017\cdots$\\ \hline
$101$&$-0.5248\cdots$&$-0.5247\cdots$&$-0.5247\cdots$&$-0.52467\cdots$&$0.00001\cdots$\\ \hline
$103$&$-0.5276\cdots$&$-0.5272\cdots$&$-0.5272\cdots$&$-0.52717\cdots$&$0.00003\cdots$\\ \hline
$107$&$-0.5262\cdots$&$-0.5256\cdots$&$-0.5260\cdots$&$-0.52609\cdots$&$0.00003\cdots$\\ \hline
$109$&$-0.5133\cdots$&$-0.5134\cdots$&$-0.5135\cdots$&$-0.51362\cdots$&$0.00002\cdots$\\ \hline
$113$&$-0.5420\cdots$&$-0.5414\cdots$&$-0.5416\cdots$&$-0.54164\cdots$&$0.00002\cdots$\\ \hline
$127$&$-0.5318\cdots$&$-0.5313\cdots$&$-0.5311\cdots$&$-0.53121\cdots$&$0.00591\cdots$\\ \hline
$131$&$-0.5556\cdots$&$-0.5555\cdots$&$-0.5556\cdots$&$-0.55564\cdots$&$0.00002\cdots$\\ \hline
$137$&$-0.5411\cdots$&$-0.5412\cdots$&$-0.5411\cdots$&$-0.54113\cdots$&$0.00003\cdots$\\ \hline
$139$&$-0.5348\cdots$&$-0.5346\cdots$&$-0.5346\cdots$&$-0.53471\cdots$&$0.00079\cdots$\\ \hline
$149$&$-0.5367\cdots$&$-0.5369\cdots$&$-0.5370\cdots$&$-0.53700\cdots$&$0.00000\cdots$\\ \hline
$151$&$-0.5433\cdots$&$-0.5436\cdots$&$-0.5438\cdots$&$-0.54378\cdots$&$0.00002\cdots$\\ \hline
$157$&$-0.5297\cdots$&$-0.5298\cdots$&$-0.5298\cdots$&$-0.52986\cdots$&$0.00006\cdots$\\ \hline
$163$&$-0.5407\cdots$&$-0.5406\cdots$&$-0.5407\cdots$&$-0.54078\cdots$&$0.00001\cdots$\\ \hline
$167$&$-0.5281\cdots$&$-0.5291\cdots$&$-0.5289\cdots$&$-0.52899\cdots$&$0.00000\cdots$\\ \hline
$173$&$-0.5575\cdots$&$-0.5576\cdots$&$-0.5576\cdots$&$-0.55769\cdots$&$0.00001\cdots$\\ \hline
\end{tabular}
\end{center}

\vfil\eject
\noindent {\tt Acknowledgement}. I'd like to thank Y. Hashimoto and Y. Ihara for kindly sending me
\cite{veel}, respectively \cite{I2} and helpful e-mail
correspondance. Furthermore, I like to thank P. Pollack for pointing out
the existence of Rosen's paper \cite{R} to me. Last but not least thanks are due to W. Bosma
for implementing formula (\ref{noidea}) in MAGMA (to use MAGMA for this
was suggested to me by H. Gangl). The columns headed `true' in Table 1 and 2
were produced using the data of Bosma.

\medskip\noindent {\footnotesize 
{\tt Pieter Moree}, Max-Planck-Institut f\"ur Mathematik,
Vivatsgasse 7, D-53111 Bonn, Germany.\\ e-mail: {\tt moree@mpim-bonn.mpg.de}}
\end{document}